\documentclass[12pt]{article}
\usepackage{amsfonts,amssymb}
\usepackage{amsmath}

\begin{document}

\begin{center}
{\Large \ About division quaternion algebras\\[0pt]
and division symbol algebras}

\begin{equation*}
\end{equation*}%
Diana SAVIN%
\begin{equation*}
\end{equation*}
\end{center}

\textbf{Abstract. }{\small In this paper, we find a class of division
quaternion algebras over the field }$\mathbb{Q}\left( i\right) ${\small \ and a class
of division symbol algebras over a cyclotomic field.}

\bigskip

\smallskip


\bigskip


\begin{equation*}
\end{equation*}

\textbf{1. Preliminaries}%
\bigskip\\
Let $K$ be a field with \textit{char} $K\neq 2.$ Let $A$ be a simple $K$-algebra and $%
Z\left( A\right) $ be the center of $A.$ We recall that the $K$- algebra $A$
is called central simple if $Z\left( A\right) =K.$\newline
Let $n$ be an arbitrary positive integer, $n\geq 3$ and let $\xi $ be a
primitive $n$-th root of unity. If $char(K)$ does not divide $n$ and $\xi $$%
\in $$K$, $\ $let $K^{\ast }=K\backslash \{0\},$ $a,b$ $\in K^{\ast }$ and
let $A$ be the algebra over $K$ generated by elements $x$ and $y$ where%
\begin{equation*}
x^{n}=a,y^{n}=b,yx=\xi xy.
\end{equation*}

This algebra is called a \textit{symbol algebra }(also known as a \textit{%
power norm residue algebra}) and it is denoted by $\left( \frac{a,~b}{%
K,\omega }\right) .$ J. Milnor, in [Mi; 71], calls it the symbol algebra. For $%
n=2,$ we obtain the quaternion algebra. Quaternion algebras and symbol
algebras are central simple algebras.
Quaternion algebras and symbol algebras have many applications in number theory (class field theory).
Conditions  of some algebras to be split or with division  were
intensively studied in various papers, \ as for example in the papers  [Fl, Sa; 15], [Sa, Fl, Ci; 09]
and [Sa; 14]; [Fl, Sa; 14] in which the authors found some interesting examples of
quaternion division algebras ,
respectively quaternion algebras and symbol algebras which split. In this
paper, using some of these results and some properties of cyclotomic fields
and $p-$ adic fields, we find a class of division quaternion algebras over
the field $\mathbb{Q}\left( i\right) $ (see Theorem 3.1) and a class of
division symbol algebras over a cyclotomic field (see Theorem 3.2). 
\begin{equation*}
\end{equation*}

\textbf{2. Introduction}%
\bigskip\\
In the following, we assume that $K$ is a commutative field and $A$ is a
finite dimensional algebra over $K.$ If $A$ is a central simple $K-$
algebra, then the dimension $n$ of $A$\ over $K$ is a square. The positive
integer $d=\sqrt{n}$ is called\textit{\ the degree} of the algebra.\newline
We recall some definitions and properties of the theory of associative
algebras, cyclotomic fields and $p-$ adic fields, which will be used in our
paper.\medskip

\textbf{Definition 2.1.} Let $A\neq 0$ be an algebra over the field $K.$ If
the equations $ax=b,\,ya=b,\forall a,b\in A,a\neq 0,$ have unique solutions,
then the algebra $A$ is called\textit{\ a division algebra. }If $A$ is a
finite-dimensional algebra, then $A$ is a division algebra if and only if $A$
is without zero divisors ($x\neq 0,y\neq 0\Rightarrow xy\neq 0$)\medskip .

\textbf{Definition 2.2.} Let $K\subset L$ be a fields extension and let 
\textit{\ }$A$\textit{\ } be a central simple algebra over the field $K.$ We
recall that: \newline
i) $A$ is called \textit{split } \textit{by} $K$ if $A$ is isomorphic with a
matrix algebra over $K.$\newline
ii) $A$ is called \textit{split } \textit{by} $L$ and $L$ is called\textit{\
a} \textit{\textbf{\ }splitting field } \textit{for} $A$ if $A\otimes _{K}L$
is a matrix algebra over $L.\medskip $

We will denote by $\mathbb{H}_{K}\left( \alpha ,\beta \right) ~$\textit{the
generalized quaternion algebra} over the field $K,$the algebra of the
elements of the form $a=a_{1}\cdot 1+a_{2}e_{2}+a_{3}e_{3}+a_{4}e_{4},$
where $a_{i}\in K,i\in \{1,2,3,4\}$, and the elements of the basis $%
\{1,e_{2},e_{3},e_{4}\}$ satisfy the following multiplication table: \vspace{%
3mm}

\begin{center}
\begin{tabular}{ccccc}
$\cdot \,\,$\vline & $1$ & $e_{2}$ & $e_{3}$ & $e_{4}$ \\ \hline
$1$\thinspace \vline & $1$ & $e_{2}$ & $e_{3}$ & $e_{4}$ \\ 
$e_{2}$\vline & $e_{2}$ & $-\alpha $ & $e_{4}$ & $-\alpha e_{3}$ \\ 
$e_{3}$\vline & $e_{3}$ & $-e_{4}$ & $-\beta $ & $\beta e_{2}$ \\ 
$e_{4}$\vline & $e_{4}$ & $\alpha e_{3}$ & $-\beta e_{2}$ & $-\alpha \beta $%
\end{tabular}
\end{center}

\bigskip

We denote by $\boldsymbol{n}\left( a\right) $ the norm of a generalized
quaternion $a.$ This norm has the following expression $\boldsymbol{n}\left(
a\right) =a_{1}^{2}+\alpha a_{2}^{2}+\beta a_{3}^{2}+\alpha \beta a_{4}^{2}.$%
This algebra is a division algebra if and only if for $x\in \mathbb{H}%
_{K}\left( \alpha ,\beta \right) $ we have $\boldsymbol{n}\left( x\right) =0$
if and only if $x=0.\,$ Otherwise, the algebra $\mathbb{H}_{K}\left( \alpha
,\beta \right) $ is called \textit{a split }algebra.\newline
In the books [Lam; 04], [Pi; 82], [Gi, Sz; 06] appear the following criterions to decide if a quaternion algebra or a symbol
algebra is split.
\medskip

\textbf{Proposition 2.1.} ([Lam; 04], [Pi; 82]) \textit{The quaternion algebra }$\mathbb{H}%
_{K}\left( \alpha ,\beta \right) $ \textit{\ is split algebra if and only } $%
\beta $\textit{\ is a norm from the extension }$K\subseteq K(\sqrt{\alpha }%
)\medskip $

\textbf{Proposition 2.2.} ([Gi, Sz; 06]) \textit{The quaternion algebra }$\mathbb{H}%
_{K}\left( \alpha ,\beta \right) $\textit{\ is split if and only if the
conic }$C\left( \alpha ,\beta \right) :$ $\alpha x^{2}+\beta y^{2}=z^{2}$ 
\textit{\ has a rational point over }$K($\textit{i.e. if there are }$%
x_{0},y_{0},z_{0}\in K$\textit{\ such that }$\alpha x_{0}^{2}+\beta
y_{0}^{2}=z_{0}^{2}).\medskip $

\textbf{Theorem 2.1.} ([Gi, Sz; 06]) \textit{Let }$K$\textit{\ be a field such that }$%
\zeta \in K,\,\,\zeta ^{n}=1,\zeta $\textit{\ is a primitive root, and let }$%
\alpha ,\beta \in K^{\ast }.$\textit{\ Then the following statements are
equivalent:}\newline
\textit{i) The cyclic algebra }$A=\left( \frac{\alpha ,\beta }{K,\zeta }%
\right) $\textit{\ is split.}\newline
\textit{ii) The element }$\beta $\textit{\ is a norm from the extension }$%
K\subseteq K(\sqrt[n]{\alpha }).\medskip $

\textbf{Theorem 2.2. (The Wedderburn norm criterion) }([Led; 05]). \textit{Let} $n
$ \textit{be a positive integer,} $n\geq 3$ \textit{and let} $L/K$ \textit{%
be a cyclic fields extension of order} $n.$ \textit{Let} $\sigma $ \textit{%
be a generator of the Galois group} Gal($L/K$). \textit{Then} $\left(
M,\sigma ,a\right) $ \textit{is a division algebra if} $a^{d}$ \textit{is
not a norm in} $L/K$ \textit{for} $d|n,d<n.\medskip $

\textbf{Theorem 2.3} (\textbf{Weddeburn}) ([Mil; 08], [Mi; 71]) \textit{Let }$A$\textit{\
be a central simple algebra over the field}$\,\,K.$\textit{\ Therefore there
are }$n\in \mathbb{N}^{\ast }$\textit{\ and a division algebra }$D,$\textit{%
\ }$K\subseteq D,$\textit{\ such that }$A\simeq $\textit{\ }$\mathcal{M}%
_{n}\left( D\right) .$\textit{\ The division algebra }$D$\textit{\ is unique
up to an isomorphism.\medskip }

\textbf{Theorem 2.4} ([Lan; 02]) \textit{Let } $K$ \textit{be a field,} $n$ 
\textit{be an arbitrary positive integer such that} $g.c.d.\left(
n,charK\right) =1$ \textit{and }$K$\textit{\ contains a primitive root of
order} $n$ \textit{of unity.}\newline
i) \textit{Let} $L$ \textit{be a cyclic extension of degree} $n.$ \textit{%
Then, there is} $\alpha \in K$ \textit{such that} $L=K\left( \alpha \right) $
\textit{and} $\alpha $ \textit{satisfies the equation} $x^{n}-a=0$ \textit{%
for some} $a\in K.$\newline
ii) \textit{Conversely, let} $a\in K$ \textit{\ and } $\alpha $ \textit{be a
root of the equation} $x^{n}-a.$ \textit{Then} $K\left( \alpha \right) $ 
\textit{is cyclic over} $K,$ \textit{of degree} $d,$ $d|n$ \textit{and} $%
\alpha ^{d}\in K.$ \medskip

In [Br, Pa; 74] E. Brown and J. Parry determined all imaginary bicyclic biquadratic fields $K=\mathbb{Q}\left(\sqrt{l}, \sqrt{d}\right)$ with class number $1.$ From these fields, we use in the section $3$ the imaginary biquadratic number fields $K=\mathbb{Q}\left(\sqrt{l}, \sqrt{d}\right)$ with $l=-1.$\medskip

\textbf{Theorem 2.5} ([Br, Pa; 74]) \textit{Let} $d<-1$ \textit{be a square free integer and the biquadratic field} $K=\mathbb{Q}\left(i,\sqrt{d}\right).$ \textit{Then, only values of d for which K has class number} $1$  \textit{are:} $d\in\left\{-163, -67, -43, -37, -19, -13, -11, -7, -5, -3, -2\right\}.$
\begin{equation*}
\end{equation*}

\textbf{3. Division quaternion algebras and symbol algebras, over a
quadratic field or over a cyclotomic field} 
\begin{equation*}
\end{equation*}

It is known that a quaternion algebra or a symbol algebra of degree $p$ is
either split or a division algebra (see [Lam; 04], [Led; 05]).

In the papers [Fl, Sa; 15], [Sa, Fl, Ci; 09], we found some examples of split quaternion and
symbol algebras over a quadratic field or over a cyclotomic field.\medskip\
We obtained the following results:\medskip 

\textbf{Proposition 3.1} ([Sa, Fl, Ci; 09]) \textit{Let} $p$ \textit{be a prime positive
integer,} $p\equiv 1$\textit{(\ mod }$3$) \textit{and let} $K=$\textit{\ }$%
\mathbb{Q}\left( \sqrt{3}\right) .$\textit{\ Then the quaternion algebra }$%
\mathbb{H}_{K}\left( -1,p\right) $ \textit{\ is a split algebra.\medskip }

\textbf{Proposition 3.2.} ([Sa, Fl, Ci; 09]) \textit{Let} $\epsilon $ \textit{be a
primitive root of order} $3$ \textit{of the unity. Then the algebras }$%
A=\left( \frac{\alpha ,\beta }{\mathbb{Q}\left( \varepsilon \right)
,\varepsilon }\right) ,$\textit{\ for }$\alpha ,\beta \in \{-1,1\}$ \textit{%
are split algebras}.\medskip 

Let $\epsilon $ be a primitive root of order $3$ of the unity and $\mathbb{Q}%
\left( \epsilon \right) $ be the cyclotomic field. In the paper [Fl, Sa; 15] using
the computer algebra system MAGMA, we obtained that the symbol
algebras $\left( \frac{7,11^{3}}{\mathbb{Q}\left( \epsilon \right) ,\epsilon 
}\right) ,$ $\left( \frac{7,\left( 11+\epsilon \right) ^{3}}{\mathbb{Q}%
\left( \epsilon \right) ,\epsilon }\right) ,$ $\left( \frac{7,5^{3}}{\mathbb{%
Q}\left( \epsilon \right) ,\epsilon }\right) $ are split  algebras and the class number of the
Kummer field $\mathbb{Q}\left( \epsilon, \sqrt[3]{7} \right) $ is $3.$  Moreover, in the paper [Fl, Sa; 15],
we found a class of split symbol algebras, over a cyclotomic field.\medskip 

\textbf{Proposition 3.3.} ([Fl, Sa; 15]) \textit{Let} $q$ \textit{be an odd prime
positive integer and} $\xi $ \textit{be a primitive root of order} $q$ 
\textit{of unity and let }$K=\mathbb{Q}\left( \xi \right) $\textit{\ be the
cyclotomic field.} \textit{Let} $\alpha \in K^{\ast },$ $p$ \textit{be} 
\textit{a prime rational integer,} $p\neq 3$ \textit{and let }$L=K\left( 
\sqrt[q]{\alpha }\right) $ \textit{be} \textit{the Kummer field} \textit{%
such that} $\alpha $ \textit{is a} $q$ \textit{power residue modulo} $p.$ 
\textit{Let} $h_{L}$ \textit{be the class number of} $L.$ \textit{Then, the
symbol algebras} $A=\left( \frac{\alpha ,p^{h_{L}}}{K,\xi }\right) $\textit{%
\ are split.} \medskip 

In this paper we find a class of quaternion division algebras or division
symbol algebras over a $p$-adic field, over a quadratic field or over a
cyclotomic field.

We consider the quadratic field $\mathbb{Q}\left( i\right) $ ($i^{2}=-1$)
and the cyclotomic field $\mathbb{Q}\left( \epsilon \right) ,$ where $%
\epsilon $ is a primitive root of order $3$ of the unity . Using the
computer algebra system MAGMA, we obtain:\newline
$A<a,b,c>:=$ QuaternionAlgebra$<$ RationalField() $|10,29>$; $a^{2};$ $b^{2};$ $a\ast b;$
Q :=Rationals(); Z:=RingOfIntegers(Q); Z; E :=QuadraticField($-1$);\newline
a :=RootOfUnity($2$); a; $Et<t>:=$PolynomialRing(E); E; $f:=t^{2}-10;$\newline
$K<b>:=$NumberField(f); K; $b^{2};$ NormEquation(K, $29$);\newline
\textbf{evaluate}\newline
$10$;  $29$; c; Integer Ring; $-1$\\
Quadratic Field with defining polynomial $.1^{2}+1$ over the Rational Field%
\newline
Number Field with defining polynomial $t^{2}-10$ over E\newline
$10$;  $-1$; false\newline
\smallskip \newline
respective\newline
\smallskip \newline
$Q:=$Rationals(); $E:=$CyclotomicField($3$); $a:=$RootOfUnity($3$); $a$; \newline
$Et<t>:=$PolynomialRing(E); $E;$ $f:=t^{3}-7;$ $K<b>:=$NumberField($f$);\newline
$K;$ $b^{3};$; NormEquation($K,29$); NormEquation($K,43$);NormEquation($K,13$);\newline
 NormEquation($K,19$);\newline
\textbf{evaluate}\newline
$zeta_3$; Cyclotomic Field of order $3$ and degree $2$\newline
Number Field with defining polynomial $t^{3} - 7$ over$E$; $7$\newline
true$[(-zeta_3 - 1)*b^{2}+ (-2*zeta_3 - 2)*b - 2*zeta_3 - 2]$\newline
false; false; false\newline
\smallskip\newline
Therefore, $10$ is not a quadratic residue modulo $29,$ $29\equiv1$ (mod $4$%
)and $29 $ is not a norm from the extension $\mathbb{Q}\left(i\right)
\subseteq \mathbb{Q}\left(i, \sqrt{10}\right).$ Using similar calculations
in Magma we obtain that $15$ is not a quadratic residue modulo $29, $ $29$
is not a norm from the extension $\mathbb{Q}\left(i\right) \subseteq \mathbb{%
Q}\left(i, \sqrt{15}\right)$ and $5$ is a quadratic residue modulo $29,$ $29$
is a norm from the extension $\mathbb{Q}\left(i\right) \subseteq \mathbb{Q}%
\left(i, \sqrt{15}\right).$  So, applying Proposition 2.1 it results that the
quaternion algebras $\mathbb{H}_{\mathbb{Q}\left(i\right)}\left( 10 , 29
\right) $ $\mathbb{H}_{\mathbb{Q}\left(i\right)}\left( 15 , 29 \right) $ are
division algebras and the quaternion algebra $\mathbb{H}_{\mathbb{Q}%
\left(i\right)}\left( 5 , 29 \right) $ is a split algebra.\newline
From the second example shown in Magma, it results that $29$ is a norm from
the extension $\mathbb{Q}\left(\epsilon\right) \subseteq \mathbb{Q}%
\left(\epsilon, \sqrt[3]{7}\right),$ but $43;13;19$ are not norms from the
extension $\mathbb{Q}\left(\epsilon\right) \subseteq \mathbb{Q}%
\left(\epsilon, \sqrt[3]{7}\right).$ So, applying Theorem 2.1 or Theorem
2.2, it results that the symbol algebra $\left( \frac{7,29}{\mathbb{Q}\left(
\epsilon \right) ,\epsilon }\right)$ is a split algebra, but $\left( \frac{%
7,43}{\mathbb{Q}\left( \epsilon \right) ,\epsilon }\right),$ $\left( \frac{%
7,13}{\mathbb{Q}\left( \epsilon \right) ,\epsilon }\right),$ $\left( \frac{%
7,19}{\mathbb{Q}\left( \epsilon \right) ,\epsilon }\right)$ are division
algebras. We remark that $29\equiv2$ (mod $3$), but $43;13;19\equiv1$ (mod $%
3 $).\newline
Let $\omega$ be a primitive root of order $5$ of the unity and let the
cyclotomic field $\mathbb{Q}\left(\omega\right)$. Similarly with previous
examples, using the computer algebra system MAGMA, we obtain that the symbol
algebra $\left( \frac{19,37}{\mathbb{Q}\left( \omega \right) ,\omega }%
\right) $ is a split symbol algebra, but $\left( \frac{19,11}{\mathbb{Q}%
\left( \omega \right) ,\omega }\right),$ $\left( \frac{19,31}{\mathbb{Q}%
\left( \omega \right) ,\omega }\right)$ are division symbol algebras. We
remark that $37\equiv2$ (mod $5$), but $11;31\equiv1$ (mod $5$).\newline

Considering these things, we obtain the following results.
In these results we use the notations: $\left( \cdot, \cdot \right)_{p}$  for the Hilbert symbol in the $p-$ adic field $\mathbb{Q}_{p},$ $\epsilon\left(\frac{\cdot, \cdot}{K}\right)_{v}$ for the Hasse invariant at a place $v$ of a field $K,$ 
$\left(\frac{\cdot}{p}\right)$ for the Legendre symbol in $\mathbb{Z},$ respective $\left[\frac{\cdot}{p}\right]$ for the Legendre symbol in $\mathbb{Z}\left[i\right].$\medskip

\textbf{Theorem 3.1.} \textit{Let} $p$ \textit{be} \textit{a prime positive
integer such that} $p\equiv 1$ (\textit{mod} $4$) \textit{and let the quadratic field } $\mathbb{Q}%
\left( i\right) $ ($i^{2}=-1$). \textit{Let} $\alpha $ \textit{be an integer which is not a quadratic
residue modulo} $p.$ \textit{Then the quaternion algebra} $\mathbb{H}_{%
\mathbb{Q}\left( i\right) }\left( \alpha ,p\right) $ \textit{is a division
algebra.}\medskip

\textbf{Proof.} Since $\alpha $ is not a quadratic residue modulo $p,$ it
results that $\overline{\alpha }$$\notin $$\left( \mathbb{F}_{p}^{\ast }\right) ^{2}.$
Therefore, $\mathbb{F}_{p}\left( \sqrt{\alpha }\right) /\mathbb{F}_{p}$ is a
cyclic extension of degree $2.$ From Hensel 's lemma ([1]), we know that the $p-$
adic field $\mathbb{Q}_{p}$ contains the roots of order $p-1$ of the unity.
Since $p\equiv 1$ (mod $4$), we have that $i$$\in $$\mathbb{Q}_{p},$
therefore $\mathbb{Q}\left( i\right) \subset \mathbb{Q}_{p}.$ We consider
the quaternion algebra $\mathbb{H}_{\mathbb{Q}\left( i\right) }\left( \alpha
,p\right) \otimes _{\mathbb{Q}\left( i\right) }\mathbb{Q}_{p}=\mathbb{H}_{%
\mathbb{Q}_{p}}\left( \alpha ,p\right) .$ \\
We consider the equation $\alpha x^{2} + py^{2}=z^{2}.$ We calculate the Hilbert symbol in $\mathbb{Q}_{p}:$ $\left(\alpha, p\right)_{p}.$ Since $\alpha $ is not a quadratic residue modulo $p,$ it results that $p$ does not divide $\alpha .$ Therefore
$\left(\alpha, p\right)_{p}=\left(\frac{\alpha}{p}\right)=-1.$ This implies that the equation $\alpha x^{2} + py^{2}=z^{2}$ does not have solutions in $p$-adic field $\mathbb{Q}_{p}.$  Applying
Proposition 2.2, it results that $\mathbb{H}_{\mathbb{Q}_{p}}\left( \alpha
,p\right) $ is not split, therefore $\mathbb{H}_{\mathbb{Q}_{p}}\left(
\alpha ,p\right) $ is a division algebra. This implies that $\mathbb{H}_{%
\mathbb{Q}\left( i\right) }\left( \alpha ,p\right) $ is a division algebra.
\medskip

A question which appears in the following is: what happens with the quaternion algebra $\mathbb{H}_{%
\mathbb{Q}\left( i\right) }\left( \alpha ,p\right) $ when $\alpha $ is a quadratic residue modulo $p.$ 
Using Theorem 2.5, the decomposition of a prime integer in the ring of integers of a biquadratic field (see [Mar; 95] ) and a reasoning similar to that which we used in the proof of Proposition 3.3 (see [Fl, Sa; 15]) we obtain: \medskip

\textbf{Proposition 3.4.}\textit{Let} $\alpha$$\in$ $\left\{  \pm 2, \pm 3, \pm 5, \pm 7, \pm 11, \pm 13, \pm 19, \pm 37, \pm 43, \pm 67, \pm 163\right\}$ \textit{and let} $p$ \textit{be} \textit{an odd prime  positive
integer such that} $\alpha $ \textit{is a quadratic residue modulo} $p$  \textit{and let the quadratic field } $\mathbb{Q}%
\left( i\right) $ ($i^{2}=-1$). \textit{Then the quaternion algebra} $\mathbb{H}_{%
\mathbb{Q}\left( i\right) }\left( \alpha ,p\right) $ \textit{is a split algebra.}\medskip

\textbf{Proof.} Our first remark is the fact that for every $\alpha$ from the set\\
$\left\{  \pm 2, \pm 3, \pm 5, \pm 7, \pm 11, \pm 13,
\pm 19, \pm 37, \pm 43, \pm 67, \pm 163\right\}$ there exists an odd prime  positive
integer $\alpha $ such that $\alpha $ is quadratic residue modulo $p.$ \\
Let $\mathcal{O}_{K}$ the ring of integers of the biquadratic field $K=\mathbb{Q}\left(i,\sqrt{\alpha}\right)= \mathbb{Q}\left(i,\sqrt{-\alpha}\right).$ 
From the hypothesis  follows immediately that $\mathcal{O}_{K}$ is a principal ring.
We know that, if $p\equiv 1$ (mod $4$), then $p$ splits in the ring $\mathbb{Z}\left[i\right]$ in a product 
of two primes from $\mathbb{Z}\left[i\right],$ respective, if $p\equiv 3$ (mod $4$), then $p$ is inert in the ring 
$\mathbb{Z}\left[i\right].$\\
\textbf{Case $1$}: if $p \equiv1$ (mod $4$).
We know that $\mathbb{Z}\left[i\right]$ is a principal ring. So, we have:
$$p\mathbb{Z}\left[i\right]=p_{1}\mathbb{Z}\left[i\right]p_{2}\mathbb{Z}\left[i\right],$$
where $p_{1}, p_{2}$ are prime elements from $\mathbb{Z}\left[i\right].$ Since $\alpha $ is quadratic residue modulo $p,$ it results that $\alpha $ is quadratic residue modulo $p_{1}, $ $p_{2}.$ So, we obtain the following decomposition of the ideal $p\mathcal{O}_{K}:$
$$p\mathcal{O}_{K}=P_{11}P_{12}P_{21}P_{21},$$
where $P_{i1}$ and $P_{i2},$ $i=\overline{1,2}$ are prime, principal conjugate ideals from the ring $\mathcal{O}_{K}.$ It results that $p=N_{K/\mathbb{Q}\left(i\right)} \left(P_{11}\right).$ But $P_{11}$ is a principal ideal, therefore, there exists $a\in K$ such that $p=N_{K/\mathbb{Q}\left(i\right)} \left(a\right).$ Applying Proposition 2.1 it results that the quaternion algebra $\mathbb{H}_{%
\mathbb{Q}\left( i\right) }\left( \alpha ,p\right) $ is a split algebra.\\
\textbf{Case $2$}: if $p \equiv3$ (mod $4$), we know that $p$ is inert in the ring $\mathbb{Z}\left[i\right]$ and having in view that $\alpha $ is quadratic residue modulo $p,$ we obtain that $p\mathcal{O}_{K}=P_{1}P_{2},$ where $P_{i}$ and $P_{i},$ $i=\overline{1,2}$ are prime, principal conjugate ideals from the ring $\mathcal{O}_{K}.$  Similarly with
 the case 1, we obtain that the quaternion algebra $\mathbb{H}_{%
\mathbb{Q}\left( i\right) }\left( \alpha ,p\right) $ is a split algebra.\\
\medskip

In the case when $p^{'}$ is a prime positive integer, $p^{'}\equiv 1$ (mod $4$) and $\alpha $ is an integer which is a quadratic residue modulo $p^{'}$ we obtain the following result: \medskip

\textbf{Proposition 3.5.} \textit{Let} $p^{'}$ \textit{be} \textit{a prime positive
integer such that}  $p^{'}\equiv 1$ (\textit{mod} $4$) \textit{and let the quadratic field } $\mathbb{Q}%
\left( i\right) $ ($i^{2}=-1$). \textit{Let} $\alpha $ \textit{be an integer such that}  $\alpha $ \textit{is a quadratic residue modulo} $p^{'}.$ \textit{Then the quaternion algebra} $\mathbb{H}_{%
\mathbb{Q}\left( i\right) }\left( \alpha , p^{'}\right) $ \textit{is a split algebra.}\medskip

\textbf{Proof.} We prove that the equation $\alpha x^{2} + p^{'}y^{2}=z^{2}$ has solutions over $\mathbb{Q}%
\left( i\right) .$ For this we determine the ramified primes in the quaternion algebra $\mathbb{H}_{%
\mathbb{Q}\left( i\right) }\left( \alpha , p^{'}\right). $  It is known that a such prime $p$ divides $2\alpha p^{'}$ ([Ko], [Ko; 00]).
Since $p^{'}\equiv 1$ (mod $4$)  it results that $p^{'}\mathbb{Z}%
\left[ i\right]=p^{'}_{1}\mathbb{Z}%
\left[ i\right]\cdot p^{'}_{2}\mathbb{Z}%
\left[ i\right],$ where $p^{'}_{1}\mathbb{Z}%
\left[ i\right], p^{'}_{2}\mathbb{Z}%
\left[ i\right]$$\in$ Spec($\mathbb{Z}%
\left[ i\right]$).\\
We calculate the Hasse invariant: $\epsilon\left(\frac{\alpha, p^{'}}{\mathbb{Q}%
\left( i\right)}\right)_{p^{'}_{j}} = \left[\frac{\alpha}{p^{'}_{j}}\right]=1,$ $j=\overline{1,2}.$ It results that 
$p^{'}_{1},$ $p^{'}_{2}$ are not ramify in $\mathbb{H}_{%
\mathbb{Q}\left( i\right) }\left( \alpha , p^{'}\right). $\\
\textbf{Case 1:} if $2$$\nmid$$\alpha .$\\
Let  $\alpha=q^{\beta_{1}}_{1}\cdot q^{\beta_{2}}_{2}\cdot...\cdot q^{\beta_{r}}_{r},$  where $r$ is an odd natural number,
$q_{j},$  $j=\overline{1,r}$ are odd prime integers and $\beta_{j}$$\in$$\mathbb{N}^{*},$ for $j=\overline{1,r}.$\\
Let $p$$\in$$\left\{q_{1}, q_{2},..., q_{r}\right\}.$ It results that $p \equiv 3$ (mod $4$) or $p \equiv 1$ (mod $4$).\\
If $p \equiv 3$ (mod $4$),  it results that $p$ remain prime in the ring $\mathbb{Z}[i].$ Using the properties of the Hasse invariant we obtain: 
$$\epsilon\left(\frac{\alpha, p^{'}}{\mathbb{Q}%
\left( i\right)}\right)_{p} =\prod\limits_{k=1, k\neq j}^{r} \epsilon\left(\frac{q^{\beta_{k}}_{k}, p^{'}}{\mathbb{Q}%
\left( i\right)}\right)_{p} \cdot \epsilon\left(\frac{p^{\beta_{j}}, p^{'}}{\mathbb{Q}%
\left( i\right)}\right)_{p}=$$
$$=\left(\epsilon\left(\frac{p, p^{'}}{\mathbb{Q}%
\left( i\right)}\right)_{p}\right)^{\beta_{j}}=\left[\frac{p^{'}}{p}\right]^{\beta_{j}}.\eqno(3.1)$$ 
Since $p^{'}$$\equiv$$1$ (mod $4$),  it results that  $p^{'}\mathbb{Z}%
\left[ i\right]=p^{'}_{1}\mathbb{Z}%
\left[ i\right]\cdot p^{'}_{2}\mathbb{Z}%
\left[ i\right],$ where $p^{'}_{1}\mathbb{Z}%
\left[ i\right], p^{'}_{2}\mathbb{Z}%
\left[ i\right]$ are prime ideals in $\mathbb{Z}%
\left[ i\right]$ and $p^{'}_{2}=\overline{p^{'}_{1}}.$\\
Taking into account (3.1) we obtain:
$$\epsilon\left(\frac{\alpha, p^{'}}{\mathbb{Q}%
\left( i\right)}\right)_{p} =\left[\frac{p^{'}_{1}}{p}\right]^{\beta_{j}}\cdot \left[\frac{\overline{p^{'}_{1}}}{p}\right]^{\beta_{j}}=1.$$
So, each divisor $p\equiv 3$ (mod $4$) of  $\alpha$  does not ramify in $\mathbb{H}_{%
\mathbb{Q}\left( i\right) }\left( \alpha , p^{'}\right). $\\
If $p \equiv 1$ (mod $4$),  it results that  $p\mathbb{Z}%
\left[ i\right]=p_{1}\mathbb{Z}%
\left[ i\right]\cdot p_{2}\mathbb{Z}%
\left[ i\right],$ where $p_{1}\mathbb{Z}%
\left[ i\right], p_{2}\mathbb{Z}%
\left[ i\right]$$\in$ Spec($\mathbb{Z}%
\left[ i\right]$). Analogously to the previous considerations, we obtain
the Hasse invariant $\epsilon\left(\frac{\alpha, p^{'}}{\mathbb{Q}%
\left( i\right)}\right)_{p_{1}}=\epsilon\left(\frac{\alpha, p^{'}}{\mathbb{Q}%
\left( i\right)}\right)_{p_{2}}=1.$ So,  $p_{1}, p_{2}$ do not ramify in $\mathbb{H}_{%
\mathbb{Q}\left( i\right) }\left( \alpha , p^{'}\right). $\\
\textbf{Case 2:} if $2$$\mid$$\alpha .$\\
We know $2=-i\left(1+i\right)^{2},$  $1+i$ is a prime element in  $\mathbb{Z}%
\left[ i\right],$ $i$$\in$U$\left(\mathbb{Z}%
\left[ i\right]\right).$  Considering the results obtained in case 1 and that 
$\prod\limits_{p}\epsilon\left(\frac{\alpha, p^{'}}{\mathbb{Q}%
\left( i\right)}\right)_{p} =1, $ is results that  $\epsilon\left(\frac{\alpha, p^{'}}{\mathbb{Q}%
\left( i\right)}\right)_{1+i} =1.$ So, $1+i$ does not ramify in $\mathbb{H}_{%
\mathbb{Q}\left( i\right) }\left( \alpha , p^{'}\right). $\\ 
From the previously proved, applying Minkovski-Hasse theorem we get that the equation $\alpha x^{2} + p^{'}y^{2}=z^{2}$ has solutions over $\mathbb{Q}%
\left( i\right),$ so applying Proposition 2.2 it results that the $\mathbb{H}_{\mathbb{Q}\left( i\right) }\left( \alpha , p^{'}\right) $ is a split algebra.\bigskip

We asked ourselves if the quaternion algebras from the statement of Proposition 3.5 split over $\mathbb{Q}.$ When a $F$- quaternion algebra splits over a quadratic field $F\left(\sqrt{w}\right),$ in the paper [Ri, Lam; 74] are given sufficient conditions for that the $F$- quaternion algebra splits over $F.$ But this conditions are given only when $w$ is totally positive; 
this is not our situation (when $F\left(\sqrt{w}\right)=\mathbb{Q}%
\left( i\right)$ ) \\
When $\alpha$ is also prime, in [Al, Ba; 04] is realized a classification of quaternion algebras $\mathbb{H}_{%
\mathbb{Q}}\left( \alpha , p^{'}\right)$ (in split algebras, respectively division algebras ) after congruences satisfied by  into 
$\alpha$ and $p^{'}$.\\
Making some computation in Magma we obtain that the answer at our question is negative. For example, if $\alpha=33$ $p=29,$ we have that $29 \equiv 1$ (mod $4$) and $33$ is a quadratic residue modulo $29.$ Using Magma we obtain that the discriminant (in fact the generator of this discriminant) of the quaternion algebra $\mathbb{H}_{\mathbb{Q}\left( i\right) }\left( 33 , 29\right) $ is $1,$ so the algebra $\mathbb{H}_{\mathbb{Q}\left( i\right) }\left( 33 , 29\right) $ splits, but the discriminant of the quaternion algebra $\mathbb{H}_{\mathbb{Q}}\left( 33 , 29\right) $ is $33,$  so the $\mathbb{H}_{\mathbb{Q}}\left( 33 , 29\right) $ is a division algebra. All ramified primes in the algebra $\mathbb{H}_{\mathbb{Q}}\left( 33 , 29\right) $ are $3$ and $11$ and we remark that $3$ and $11$ are not quadratic residues modulo $29.$ In another example: if $\alpha=35,$ $p=29$ we obtain that the quaternion algebra $\mathbb{H}_{\mathbb{Q}\left( i\right) }\left( 35 , 29\right) $ splits, also the quaternion algebra $\mathbb{H}_{\mathbb{Q}}\left( 35 , 29\right) $ splits; $35$ is a quadratic residue modulo $29$ and  $5;$ $7$ are also quadratic residues modulo $29.$ \\
Considering these things, we get the following result.\bigskip

\textbf{Proposition 3.6.} \textit{Let} $p^{'}$ \textit{be} \textit{a prime positive
integer such that}  $p^{'}\equiv 1$ (\textit{mod} $4$). \textit{Let} $\alpha $ \textit{be an integer such that each divisor of}  $\alpha $ \textit{is a quadratic residue modulo} $p^{'}.$ \textit{Then the quaternion algebra} $\mathbb{H}_{%
\mathbb{Q}}\left( \alpha , p^{'}\right) $ \textit{is a split algebra.}\medskip

\textbf{Proof.} The proof is similar to the proof of Proposition 3.5. . The only difference is when instead of the relation (3.1) from the proof of Proposition 3.5 appears the following situation (using the properties of the Hilbert symbol):
$$\left(\alpha, p^{'}\right)_{p} =\prod\limits_{k=1, k\neq j}^{r} \left(q^{\beta_{k}}_{k}, p^{'}\right)_{p}
\cdot \left(p^{\beta_{j}}, p^{'}\right)_{p}=\left(\left(p, p^{'}\right)_{p}\right)^{\beta_{j}}
=\left(\frac{p^{'}}{p}\right)^{\beta_{j}}.$$ 
Since $p^{'}$$\equiv$$1$ (mod $4$), applying quadratic reciprocity law, it results $\left(\frac{p^{'}}{p}\right)=
\left(\frac{p}{p^{'}}\right).$
Since each divisor of $\alpha $ is a quadratic residue modulo $p^{'},$ it results that $\left(\frac{p}{p^{'}}\right)=1.$ 
We obtain that $\left(\alpha, p^{'}\right)_{p}=1$ so, each divisor $p\equiv 3$ (mod $4$) of  $\alpha$  does not ramify in $\mathbb{H}_{%
\mathbb{Q}}\left( \alpha , p^{'}\right). $\bigskip

Now, we generalize the Theorem 3.1 for the symbol algebras.\medskip

\textbf{Theorem 3.2.} \textit{Let} $p$ \textit{and} $q$ \textit{be} \textit{%
prime positive integers such that} $p\equiv 1$ (\textit{mod} $q$), $\xi $ 
\textit{be a primitive root of order} $q$ \textit{of unity and let }$K=%
\mathbb{Q}\left( \xi \right) $\textit{\ be the cyclotomic field.} \textit{%
Then there exists an integer} $\alpha $ \textit{not divisible by} $p$ 
\textit{whose residue class mod} $p$ \textit{does not belongs to} $\left( 
\mathbb{F}_{p}^{\ast }\right) ^{q}$ \textit{and for every such an} $\alpha $%
, \textit{we have:}\newline
i) \textit{the algebra} $A\otimes _{K}\mathbb{Q}_{p}$ \textit{is a division
algebra over} $\mathbb{Q}_{p},$ \textit{where} $A$ \textit{is the symbol
algebra} $A=\left( \frac{\alpha ,p}{K,\xi }\right) ;$\newline
ii) \textit{the symbol algebra} $A$ \textit{is a division algebra over} $%
K.\medskip $

\textbf{Proof.} Let be the homomorphism $f:\mathbb{F}_{p}^{\ast }\rightarrow 
\mathbb{F}_{p}^{\ast },$ $f\left( x\right) =x^{q}.$ Since $q$ divides $p-1,$
it results $Ker\left( f\right) =\left\{ x\in {F}_{p}^{\ast }|x^{q}\equiv 1\
(mod\ p)\right\} $ is non -trivial, so $f$ is not injective. So, $f$ is not
surjective. It results that there exists $\overline{\alpha }$ (in $\mathbb{F}_{p}^{\ast
}, $) which does not belongs to $\left( \mathbb{F}_{p}^{\ast }\right) ^{q}.$ So, $\mathbb{F}_{p}\subset\mathbb{F}%
_{p}\left( \sqrt[q]{\overline{\alpha }}\right).$
\newline
 The extension of fields $\mathbb{F}%
_{p}\left( \sqrt[q]{\overline{\alpha }}\right) /\mathbb{F}_{p}$ is a cyclic
extension of degree $q.$ Applying Hensel 's lemma and the
fact that $p\equiv 1$ (mod $q$), it results that $\mathbb{Q}_{p}$ contains
the $q$-th roots of the unity, therefore $\mathbb{Q}\left( \xi \right)
\subset \mathbb{Q}_{p}.$ We consider the symbol algebra $A\otimes _{K}%
\mathbb{Q}_{p}=\left( \frac{\alpha ,p}{\mathbb{Q}_{p},\xi }\right) .$  Applying Theorem 2.4, it results that the
extension $\mathbb{Q}_{p}\left( \sqrt[q]{\alpha }\right) /\mathbb{Q}_{p}$ is
a cyclic unramified extension of degree $q,$ therefore a norm of an element
from this extension can be a positive power of $p,$ but can not be $p.$
Applying Theorem 2.1, we get that $\left( \frac{\alpha ,p}{\mathbb{Q}%
_{p},\xi }\right) $ is not a split algebra, therefore it is a division
algebra. This implies that $A$ is a symbol division algebra.\bigskip

\textbf{Conclusions.} In this paper we found a class of quaternion division
algebras or division symbol algebras over a $p$-adic field, over a quadratic
field or over a cyclotomic field.

Using the computer algebra system MAGMA over the quadratic field $\mathbb{Q}%
\left( i\right) $ ($i^{2}=-1$) and the cyclotomic field $\mathbb{Q}\left(
\epsilon \right) ,$ where $\epsilon $ is a primitive root of order $3$ of
the unity, we obtain very good examples which allowed us to find conditions
in the Theorem 3.1, Proposition 3.4, Proposition 3.5 and Theorem 3.2. \bigskip

\textbf{Acknowledgements.} The author is very grateful to Professor Victor
Alexandru for many helpful discussions about this paper which helped the author to improve this paper.  The author thanks 
Professor David Kohel for the discussions about the ramified primes in a quaternion algebra and Professor Ali Mouhib for the discussions about biquadratic fields. Also, the author thanks Professors Ezra Brown and Kenneth S. Williams for the fact that they provided to the author the paper [Br, Pa; 74] and Professor Montse Vela for the fact that she provided to the author the paper [Ri, Lam; 74].
\begin{equation*}
\end{equation*}
\textbf{References}%
\begin{equation*}
\end{equation*}%
[Al, Go; 99] V. Alexandru, N.M Gosoniu, \textit{Elements of Number Theory } (in Romanian), Ed.
 Bucharest University, 1999.\newline
[Al, Ba; 04] M. Alsina, P. Bayer, \textit{Quaternion Orders, Quadratic Forms and Shimura Curves,} CRM Monograph Series, vol. \textbf{22},  American Mathematical Society, 2004. \newline
[Br, Pa; 74] E. Brown, C. J. Parry, \textit{The imaginary bicyclic biquadratic fields with class - number 1,} J. Reine Angew Math. \textbf{226}, 1974, p. 118-126.\newline
[Ri, Lam; 74]  R. Elman,  T.Y. Lam, \textit{Classification Theorems for Quadratic Forms over Fields,} Commentarii Mathematici Helvetici, \textbf{49}, 1974, p. 373-381.\newline
[Fl, Sa; 14] C. Flaut, D. Savin, \textit{Some properties of the symbol algebras of
degree} $3$, Math. Reports, vol. \textbf{16(66)}•, no. 3, • 2014, p.443-463.\newline
[Fl, Sa; 15] C. Flaut, D. Savin, \textit{Some examples of division symbol algebras of
degree }$\mathit{3}$ \textit{and} $5$, accepted in Carpathian J Math.\newline
[Gi, Sz; 06] P. Gille,  T. Szamuely, \textit{Central Simple Algebras and Galois
Cohomology}, Cambridge University Press, 2006.\newline
[Ko] D. Kohel,  \textit{Quaternion algebras,} echidna.maths.usyd.edu.au/kohel/alg/doc/\newline
AlgQuat.pdf\newline
[Ko; 00] D. Kohel,  \textit{Hecke module structure of quaternions,} Proceedings of Class Field Theory - Centenary and Prospect (Tokyo, 1998), K. Miyake, ed., Advanced Studies in Pure Mathematics, \textbf{30}, 177-196, 2000. \newline
[Ko, La, Pe, Ti; 14] D. Kohel,  K. Lauter, C. Petit, and J.-P. Tignol, \textit{ On the quaternion} $l$-\textit{isogeny path problem,} LMS Journal of Computational Mathematics, \textbf{17}, 418-432, 2014. \newline
[Lam; 04] T. Y. Lam, \textit{Introduction to Quadratic Forms over Fields,}
American Mathematical Society, 2004.\newline
[Lan; 02] S. Lang, \textit{Algebra,} Springer-Verlag, 2002.\newline
[Led; 05] A. Ledet, \textit{Brauer Type Embedding Problems }, American
Mathematical Society, 2005.\newline
[Lem;00] F. Lemmermeyer, \textit{Reciprocity laws, from Euler to Eisenstein },
Springer-Verlag, Heidelberg, 2000.\newline
[Mar; 95] D. Marcus, \textit{Number fields}, Universitext, 1995. \newline
[Mil; 08] J.S. Milne, \textit{Class Field Theory},
http://www.math.lsa.umich.edu/~jmilne.\newline
[Mi; 71] J. Milnor, \ \textit{Introduction to Algebraic K-Theory, }Annals of
Mathematics Studies, Princeton Univ. Press, 1971.\newline
[Pi; 82] Pierce, R.S., \textit{Associative Algebras}, Springer Verlag, 1982.%
\newline
[Sa, Fl, Ci; 09] D. Savin, C.Flaut, C.Ciobanu, \textit{Some properties of the symbol
algebras}, Carpathian Journal of Mathematics , \textbf{25(2)(2009)}, p.
239-245.\newline
[Sa; 14] D. Savin, \textit{About some split central simple algebras}, An. Stiin. Univ. "`Ovidius" Constanta, Ser.
Mat, \textbf{22} (1) (2014), p. 263-272.
\[\]

{\footnotesize \vspace{2mm} \noindent 
\begin{minipage}[b]{10cm}
Diana SAVIN, \\
Faculty of Mathematics and Computer Science, \\ 
Ovidius University of Constanta, \\ 
Constanta 900527, Bd. Mamaia no.124, Rom\^{a}nia \\
Email: savin.diana@univ-ovidius.ro\end{minipage}
}

{\footnotesize \bigskip }

\end{document}